\input amstex

\documentstyle{amsppt}
\NoBlackBoxes


\define\ga{\gamma}

\define\bck{{\backslash}}
\define\ind{\operatorname{ind}}

\define\maj{\operatorname{maj}}

 at 8truept   
 at 10truept   
 at 10truept
 at 8truept

\define\today{\ifcase\month\or January\or February\or March\or
  April\or May\or June\or July\or August\or September\or
  October\or November\or December\fi
  \space\number\day, \number\year}

\rightheadtext{Congruence classes of major index}
\topmatter
\title On counting permutations by pairs of congruence classes of major 
index
\endtitle
\author H\'el\`ene Barcelo, Robert Maule and Sheila Sundaram\endauthor

\address Department of Mathematics, Arizona State University, 
Tempe, AZ\endaddress
\email barcelo\@asu.edu\endemail

\address Department of Mathematics, Arizona State University, 
Tempe, AZ\endaddress
\email rgmaule\@msn.com\endemail

\address 240 Franklin Street Extension, Danbury, CT 06811\endaddress
\email sheila\@claude.math.wesleyan.edu\endemail



\date \today\enddate



\abstract 

For a fixed positive integer $n,$ let $S_n$ denote the symmetric group 
of $n!$ permutations on $n$ symbols, and let ${\maj(\sigma)}$  
denote the major index of a permutation $\sigma.$  
Fix positive integers $k<\ell\leq n,$ and nonnegative integers $i,j.$
Let $m_n(i\bck k; j\bck \ell)$ denote the cardinality of the set 
$\{\sigma\in S_n: {\maj}(\sigma)\equiv i \mod k, 
                             {\maj}(\sigma^{-1})\equiv j \mod \ell\}.$
In this paper we give some enumerative formulas for these numbers.
  When $\ell$ divides $(n-1)$ and $k$ divides $n,$
we show that for all $i,j,$ 
$$ m_n(i\bck k; j\bck \ell)=  {n!\over k\cdot\ell} .$$

\endabstract
\endtopmatter

\document

\head{ 1. Introduction } \endhead
  Denote by $S_n$ the symmetric group of all $n!$ permutations on the 
$n$ symbols $1,\ldots , n.$   First recall some combinatorial definitions 
pertaining to permutations.  See, e.g.,  [4]. 

\proclaim{Definition 1.1}
 Let $\sigma\in S_n.$ For  $1\leq i\leq n-1,$  $i$ 
 is said to be a descent of $\sigma$ if 
 $\sigma(i)>\sigma(i+1).$  
\endproclaim

\proclaim{Definition 1.2}
The major index of $\sigma,$ denoted 
 $\maj(\sigma),$ is the sum of the descents of $\sigma.$ 
\endproclaim

 The values of the statistic 
 $\maj$ range from 0 (for the identity) to ${n\choose 2}.$   

In \cite{1}, the following result was discovered 
using certain representations of the symmetric group $S_n,$  and 
then proved by means of a bijection as well.

\proclaim{Proposition 1.3}(\cite{1, Theorem 2.6}) 
Fix an integer $0\leq i\leq n-1.$  
$$(n-1)!=|\{\sigma\in S_n: \maj(\sigma)\equiv i\mod n\}|.$$
\endproclaim 

This paper is similarly motivated by the algebraic discovery 
(\cite{S. Sundaram, unpublished}) of the identity 
$$ (n-2)! = | \{\sigma\in S_n: {\maj}(\sigma)\equiv i \mod n, 
                             {\maj}(\sigma^{-1})\equiv j \mod (n-1) \}|,
\tag A$$
\noindent where $i,j$ are fixed nonnegative integers.  

The paper is organised as follows.  In Section 2  the main 
technical lemmas are presented.  In Section 3 we derive the 
enumerative formulas, and in Section 4 we give purely bijective proofs 
of special cases of Theorem 3.1 and Proposition 2.5.

\head{ 2. Preliminaries } \endhead

This section contains the main lemmas that are needed for the rest of 
the paper.

Let $\gamma\in S_n$ be the $n$-cycle which takes $i$ to $i+1$ modulo 
 $n,$ for all $i.$  We will sometimes write $\ga_n$ for clarity. 
 The circular class of $\sigma$ is the set of 
 permutations $[\sigma]=\{\sigma\gamma^i, 0\leq i \leq n-1\}.$
The following observation is due to Klaychko \cite{3}.  For our purposes 
it is more convenient to state the result in terms of the inverse 
permutation.  This formulation also admits an easy proof, which we 
give below for the sake of completeness.

\proclaim{Lemma 2.1}(\cite{3}, \cite{2, Lemma 4.1}) Let $\sigma\in S_n.$  
Then the function $\tau\mapsto  
  {\maj}(\tau^{-1})$ takes on all 
  $n$ possible values modulo $n$ in the circular class of $\sigma.$ 
  More precisely, we have that ${\maj}(\sigma\gamma^i)^{-1}=
  {\maj}(\sigma^{-1})+i \mod n, 0\leq i \leq n-1.$
\endproclaim

\demo{Proof} Let  $\tau=a_1\ldots a_n$ (written as a word).
  Then $\tau\gamma=a_2\ldots a_n a_1.$   
  Note that $i$ is a descent of  
 $\tau^{-1}$ if and only if $i$ appears to the right of $i+1$ when 
 $\tau$ is written as a word in $\{1,2,\ldots , n\}.$ 

  By looking at occurrences of $i$ to the right of $i+1,$ it is easy to 
  see that 
$ \maj(\tau\gamma)^{-1}-\maj(\tau^{-1})=a_1-(a_1-1)=+1, $
  if $a_1\ne 1, a_1\neq n. $  
  If $a_1=1,$ then clearly $\maj(\tau\gamma)^{-1}-\maj(\tau^{-1})=+1, $
  while if $a_1=n,$ then the difference is $-(n-1).$

  Hence in all cases the difference is $+1$ modulo $n.$ \qed
\enddemo

\proclaim{Lemma 2.2} Let $\sigma\in S_{n-1},$ and let $\sigma_i$ denote 
  the permutation in $S_n$ obtained by inserting $n$ in position $i$ 
  of $\sigma, $ $1\le i\le n.$ 
   Then  for each $k$ 
   between 1 and $n,$ the values of the major index on the set 
   $\{\sigma_1,\ldots, \sigma_k\}$ form a consecutive segment of 
   integers $[m+1, m+k],$ and the value of $\maj(\sigma_{k+1})$ is either     
    $m$ or  $m+k+1$ according as $k$ is a descent of $\sigma$ or not, 
  respectively. Note that $\maj(\sigma_n)=\maj(\sigma).$ 

  In particular, on the set $\{\sigma_1,\ldots, \sigma_n\},$ 
the function $\maj$ takes on each of the $n$ values in the  
  interval $[\maj(\sigma), \maj(\sigma)+(n-1)].$ 
\endproclaim          
   
\demo{Proof} 
\comment
 If $\sigma $ is the identity, the statement is obvious, 
so we assume 
$\sigma$ has at least one descent.  
Let $\sigma=a_1\ldots a_{n-1},$ with descents in positions 
$i_1,\ldots i_d.$  Then $\sigma_k=a_1\ldots 
a_{k-1}\;  n\; a_k a_{k+1}\ldots a_{n-1}$  
for $k=2,\ldots , n-1,$ and $\sigma_1=n\; a_1\ldots a_{n-1},$ 
$\sigma_n=a_1,\ldots a_{n-1}\; n.$

First consider the  following example in order to illustrate 
the idea of the proof. 
Let $\sigma=15234,$ and consider inserting 6 in positions  
$1,\ldots, 6,$ in succession.  By direct computation, 
$$\maj(\sigma_1)=4, \{\maj(\sigma_1), \maj(\sigma_2)\}=\{4,5\},
\{\maj(\sigma_1),\ldots, \maj(\sigma_3)\}=\{3,4,5\},$$
$$\{\maj(\sigma_1),\ldots, \maj(\sigma_4)\}=\{3,4,5,6\},$$
$$\{\maj(\sigma_1),\ldots, \maj(\sigma_5)\}=\{3,4,5,6,7\},$$
\noindent  and finally 
$\{\maj(\sigma_1),\ldots, \maj(\sigma_6)\}=\{2,3,4,5,6,7\}.$

Our induction hypothesis is that on the set of $k$ 
 permutations $\sigma_i, i=1,\dots, k,$ $\maj$ takes the $k$ values 
in the interval $[m_k+1, m_k+k],$ and further that if $a_{k-1}<a_k,$ then 
the value taken on the last permutation, namely $\maj(\sigma_k),$ 
is $m_k+k,$ while if  $a_{k-1}>a_k,$ then $\maj(\sigma_k)=m_k+1.$

Note that this holds vacuously for $k=1.$  Consider now what happens for 
$k=2.$ Clearly $\maj(\sigma_1)=\maj(\sigma)+1+d\ge 3,$ since each of 
the $d$ descents in $\sigma$ is shifted to the right by 1.
Let $\Delta_k=\maj(\sigma_{k+1})-\maj(\sigma_{k}), 1\le k\le n-1.$
It is easy to see that $\Delta_1=-1$ if $a_1>a_2,$ and $\Delta_1=+1$ if 
$a_1<a_2.$
Hence $\maj$ takes all  values in the interval 
$[\maj(\sigma_1), \maj(\sigma_1)+1]$ 
or in the interval $[\maj(\sigma_1)-1, \maj(\sigma_1)]$ on the permutations 
$\{\sigma_1, \sigma_2\},$ according as $a_1<a_2$ or $a_1>a_2.$  

The induction hypothesis therefore holds for $k=2.$ 
Assume it holds for all values $j\le k.$ 

We have $$\sigma_k=a_1\ldots a_{k-1} \ n \ a_k a_{k+1}\ldots a_{n-1},$$
\noindent and 
        $$\sigma_{k+1}=a_1\ldots a_{k-1} a_k\ n\  a_{k+1}\ldots a_{n-1}.$$
As the symbol $n$ moves from position $k$ to position $(k+1)$, its
 contribution to the major index changes by $-k+(k+1).$  
There are also possible contributions from $a_{k-1}$ and $a_k.$ It   
is routine to verify further that 
$\Delta_k= +1$ if $a_{k-1}<a_k$ and $a_k<a_{k+1};$
$\Delta_k= -k$ if $a_{k-1}<a_k$ and $a_k>a_{k+1};$
$\Delta_k= +k$ if $a_{k-1}>a_k$ and $a_k<a_{k+1};$ and finally 
$\Delta_k= -1$ if $a_{k-1}>a_k$ and $a_k>a_{k+1}.$

Hence  the $(k+1)$ values of $\maj(\sigma_i), i\le k+1,$ 
are those in the interval: 
 $$[m_k+1, m_k+k]\cup\{m_k+k+1\}=[m_k+1, m_k+k+1], {\text\ if\ } 
a_{k-1}<a_k,\  a_k<a_{k+1};$$
 $$[m_k+1, m_k+k]\cup\{m_k\}=[m_k, m_k+k], {\text\ if\ } 
a_{k-1}<a_k,\  a_k>a_{k+1};$$
  $$[m_k+1, m_k+k]\cup \{m_k+1+k\}=[m_k+1,m_k+k+1], {\text\ if\ } 
a_{k-1}>a_k,\  a_k<a_{k+1};$$ 
\noindent and finally 
   $$[m_k+1, m_k+k]\cup\{m_k\}=[m_k,m_k+k], {\text\ if\ } 
a_{k-1}>a_k,\  a_k>a_{k+1}.$$
Since the value of $\maj(\sigma_{k+1})$ changes as predicted by 
the hypothesis, depending on whether $a_k<a_{k+1}$ or not, 
 by induction the result is established for all $k.$ \qed
\endcomment

Let $\sigma'=a_1\ldots a_{n-1},$ with descents in positions 
$i_1,\ldots i_d.$  Let $\sigma$ be the permutation in $S_n$ obtained by 
appending $n$ to $\sigma'.$ Hence $\maj(\sigma)=\maj(\sigma').$ 
We shall show that the value of $\maj(\sigma')$ 
increases successively by 1 as $n$ is inserted into $\sigma'$ in the
 following order:
\item{(1)} first in the positions immediately following a descent, 
starting with the right-most descent and moving  to the left;
\item{(2)} then in the remaining positions, beginning with position 1, 
from left to right. 

For instance, if $\sigma'=14253,$ then the resulting permutations, 
beginning with $\sigma$ and then in 
the order specified above, are 
$$142536, 142563, 146253, 614253, 164253, 142653,$$ 
\noindent with respective major indices 6,7,8,9,10,11.

Let $\sigma_k$ denote the permutation in $S_n$ obtained from 
$\sigma'$ by inserting $n$ in position $k.$  Thus 
 $\sigma_k=a_1\ldots 
a_{k-1}\;  n\; a_k a_{k+1}\ldots a_{n-1}$  
for $k=2,\ldots , n-1,$ and $\sigma_1=n\; a_1\ldots a_{n-1},$ 
$\sigma_n=a_1,\ldots a_{n-1}\; n.$  Let $\Delta_k$ denote the 
difference $\maj(\sigma_k)-\maj(\sigma').$

The following facts are easily verified:
\item{(1)} If $n$ is inserted immediately after a descent of $\sigma',$ i.e., 
if $k=i_j+1,\ 1\leq j\leq d,$ then $n$ contributes a descent in position 
$i_j+1,$ but the $i_j$th element ceases to be a descent.  Also the $(d-j)$ 
descents to the right of $n$ are shifted further to the right by one.
Thus
$$\Delta_k=(d-j)+(i_j+1) -i_j=d-j+1,$$
\noindent and hence the difference $\Delta_k$ ranges from 1 through $d.$ 
\item{(2)} If $1\leq k\leq i_1,$ then the $d$ descents to the right are 
shifted over by 1, and thus 
$$\Delta_k=d+k,$$
\noindent and hence $\Delta_k$ ranges from $d+1$ through $d+i_1.$ 
\item{(3)} If $n$ is inserted in position $k$ between two descents, but 
not immediately following a descent, i.e., 
if $1+i_j<k\leq i_{j+1},\ j\leq d-1,$ then 
$$\Delta_k=(d-j)+k,$$
\noindent and hence $\Delta_k$ ranges from $(d-j+2)+i_j$ through $
d-j+i_{j+1}.$ 
\item{(4)} Finally when $i_d+2\leq k\leq n-1,$ 
$$\Delta_k=k,$$
\noindent and hence $\Delta_k$ ranges from $i_d+2$ through $n-1.$ 

This establishes the claim.  It also shows that as $n$ is inserted into 
$\sigma'$ from left to right, the difference in major index goes up 
(from $\maj(\sigma')$) first by $(d+1),$ then up by one at each step,  
except when 
it is inserted immediately after the $jth$ descent, in which case it goes 
down to $(d-j+1).$  Since when $n$ is in position $n,$ $\maj(\sigma')$ 
is unchanged, this establishes the statement of the lemma.\qed
\enddemo

\proclaim{Remark 2.3} Note that in Lemma 2.2, it need not be 
 true that the values 
of $\maj$ on an arbitrary set $\{\sigma_j,\ldots,\sigma_{j+r}\},$ $j>1,$ 
form a consecutive set of integers. 
\endproclaim

\proclaim{Lemma 2.4} Let $\sigma\in S_{n-1},$ and let $\sigma_i$ denote the 
permutation in $S_n$ obtained by inserting $n$ in position $i$ of $\sigma,$ 
for $1\leq i\leq n.$  Then $\maj(\sigma_i^{-1})=\maj(\sigma)\mod (n-1).$
\endproclaim

\demo{Proof}  Consider the effect of inserting $n$ on the set of descents 
of $\sigma^{-1}.$ If $n$ is inserted to the right of $(n-1),$ there is no 
change; if $n$ is inserted to the left of $(n-1),$ then $(n-1)$ becomes 
a descent of $\sigma_i^{-1}.$  In either case, the major index of the 
inverse permutation is unchanged modulo $(n-1).$ \qed
\enddemo

\comment
Following \cite{2}, for $\sigma\in S_n,$ define 
$$\ind(\sigma)=\sum \{i:1\leq i\leq n, \ \sigma(i)>\sigma(i+1)\ {\text or\ } 
\sigma(n)>\sigma(1) {\ {\text if}\ } i=n\}.$$ 
\noindent Thus $\ind(\sigma)$ is simply the number of descents 
of $\sigma,$ plus one if $\sigma(n)>\sigma(1).$  
Note that $\ind(\sigma\ga_n)
=\ind(\sigma)$ for all $\sigma,$ i.e., $\ind$ is constant on the class 
of circular rearrangements of $\sigma.$  From the fact that 
$\maj(\sigma)=\sum \{i:1\leq i\leq n-1,\ \sigma(i)>\sigma(i+1)\},$ it 
 follows by direct computation that  

\proclaim{Lemma 2.5}(\cite{2, Lemma 4.1})  Let $\sigma\in S_n.$  Then 
$$\maj(\sigma\ga_n)=\maj(\sigma)-\ind(\sigma) \mod n.$$
\endproclaim

\endcomment

Finally  we shall need the following result, which generalises 
Proposition 1.3. It is perhaps known, although we do not know of a precise 
reference.  There is an easy generating function proof which we 
include for the sake of completeness.  In Section 4 we will give a 
constructive proof of the equivalent statement for inverse 
permutations.  

\proclaim{Proposition 2.5} 
$${n!\over k}=|\{\sigma\in S_n: \maj(\sigma)\equiv j\mod k\}|.$$
\endproclaim

\demo{Proof} Recall the well-known formula (see \cite{4})
$$\sum_{\sigma\in S_n} q^{\maj(\sigma)}=\prod_{i=1}^{n-1}(1+q+\ldots+q^i)
\tag B$$
\noindent Note that Lemma 2.2 gives an immediate inductive 
proof of formula (B).

 Now fix integers $1\le k\le n$ and $0\le j \le k-1.$ 
 To show that the number of permutations in $S_n$ 
with major index congruent to $j\mod k$ is $n!/k,$
it suffices to show that, modulo the polynomial $(1-q^k),$  the left-hand 
side of (B) equals $(n!/k) \cdot (1+q+\ldots q^{k-1}).$ 

Since $1+q+\ldots +q^i=(1-q^{i+1})/(1-q),$ it follows from the 
generating function that for 
fixed $k\le n,$ the sum on the left-hand side vanishes at all 
$k$th roots of unity not equal to 1.  Hence, modulo $(1-q^k),$ there is 
a constant $c$ such that 
$$\sum_{\sigma\in S_n} q^{\maj(\sigma)}=c(1+q+\ldots+q^{k-1}).$$
\noindent Putting $q=1$ yields $c=n!/k,$ as required. \qed
\enddemo

\head{ 3. Enumerative Results} \endhead

\comment
\proclaim{Proposition 3.5} Fix $2\leq k_1\leq n-k_2\leq n-2.$ 
Then $(n-2)!$ is the number of permutations $\sigma\in S_n$ with 
$$\maj(\sigma)\equiv i_1\mod k_1, \ \ \ 
\maj(\sigma^{-1})\equiv i_2\mod k_2, $$
\noindent and 
$$1\leq \sigma^{-1}(n)\leq k_1, \ \ 
n-k_2+1\leq \sigma^{-1}(n-1)\leq n .$$

\endproclaim

\demo{Proof}  For this proof, a better picture of the bijection of 
of Theorem 3.1 is obtained by describing it in terms of an 
 equivalence relation on $S_n.$ 

Two permutations $\alpha, \beta \in S_n$ are equivalent if, 
when $n$ is erased, they belong to the same circular class  in $S_{n-1}.$ 
Let $\ga_{n-1}$ denote the $(n-1)$-cycle $(1,2,\ldots,n-1).$ 
The equivalence classes can thus be indexed by the set $A_{n-1}$ of 
permutations $\sigma\in S_{n-1}$  
such that $\sigma(n-1)=n-1.$  It is convenient to view the class 
of one such $\sigma$ as a rectangular array of 
permutations $\{\tau_{i,j}\}_{1\leq i\leq n, 1\leq j\leq n-1}$ in $S_n$ 
with $(n-1)$ columns and $n$ rows.  We label columns from left to right, 
and rows from top to bottom.  See Example 3.6 below for the 
two equivalence classes in $S_4.$ The elements in the $n$th (bottom) row 
are simply the $(n-1)$ circular rearrangements 
$\sigma\ga_{n-1}^j,\ 0\leq j\leq n-2,$ viewed as 
permutations in $S_n$ which fix 
the symbol $n,$  and hence have $n$ inserted in position $n.$  
Note that (in this row) 
the permutation in column $j$ has $(n-1)$ in position $(n-j).$
To construct the $(n-1)$ elements in the $i$th row,  
simply insert $n$ in position 
$i$ in each of these $(n-1)$ circular rearrangements.  

Hence the permutation $\tau_{i,j}$ in column $j$ and row $i$ of 
the class indexed by 
$\sigma$ is obtained by taking the unique circular rearrangement of 
$\sigma$ which has $(n-1)$ in position $(n-j),$ (when viewed as a 
permutation of $S_{n-1}$), viz., $\sigma\ga_{n-1}^j,$ and then inserting 
$n$ into it in position $i.$ 

Now consider the permutation $\tau_{i,j}$ in row  $i$
and  column $j.$  We make three observations:

\item{ (1)} $n$ appears in position $i$ in $\tau_{i,j}.$
\item{(2)} If $i\leq n-j$ then $n$ appears to the left of $(n-1)$ in 
$\tau_{i,j}.$ 
\item{(3)} If $i\leq n-j,$ then in particular 
$(n-1)$ appears in column $(n-j+1)$ of $\tau_{i,j}.$ 
(Otherwise $(n-1)$ is in column  $(n-j)$ of $\tau_{i,j}.$) 

It follows by Lemma 2.1 
that $\maj((\sigma\ga_{n-1}^j)^{-1})=\maj(\sigma^{-1})+j,$ 
for $0\leq j\leq n-2.$  Also from (2) above we have that 
$\maj(\tau_{i,j}^{-1})=\maj((\sigma\ga_{n-1}^j)^{-1})+n-1$ for $i\leq n-j.$
Hence for fixed $i\leq n-j,$ the permutations in the first $j$ columns 
and $i$th row of the rectangular array have the property 
 that the major index of their inverses form a consecutive 
segment of integers, and in particular, form a complete residue 
system modulo $j.$ 

Now let $k_1, k_2,$ $i_1, i_2,$ be as in the statement.  It follows 
from this description of the rectangular array, and the above 
discussion, 
that, for rows $i\leq n-k_2,$  the first 
$k_2$ columns form a complete residue system modulo $k_2$ for the 
major index of the inverse permutation.  Also from Lemma 2.2 it 
follows that in each column, the first $k_1$ rows form a complete residue 
system modulo $k_1,$ for the major index of the permutation itself.

It is now easy to see that for each $\sigma\in A_{n-1},$ there is a unique 
permutation in the first $k_2$ columns and the first $k_1$ rows, 
provided $k_1\leq n-k_2,$
which satisfies the conditions in the statement of the proposition.  
It is easy to see how to recover $\sigma$ from such a permutation, 
and hence we have a bijection establishing the identity. \qed

\enddemo

\proclaim{Example 3.6}  Here are the two equivalence classes in $S_4,$ 
under the relation used in the proof of Theorem 3.1, and described explicitly 
in the proof of Proposition 3.5.  Next to each permutation $\tau$ is the  
ordered pair $(\maj(\tau), \maj(\tau^{-1})).$  Note that, as predicted 
by Lemma 2.2, in any one column, the first coordinate of the ordered 
pairs forms a consecutive segment of integers.   Also, as established 
in the proof of Proposition 3.5, in any row $i,$ the second coordinates of 
the ordered pairs in columns $1, \ldots, (n-i)$ form a consecutive 
segment of integers.

\endproclaim

$$\spreadmatrixlines{.5\jot} 
\matrix  &4123 \ (1,3) &4231\ (4,4)  &4312\ (3,5) \\
         &1423 \ (2,3) &2431\ (5,4)  &3412\ (2,2) \\
         &1243 \ (3,3) &2341\ (3,1)  &3142\ (4,2) \\ 
         &1234 \ (0,0) &2314\ (2,1)  &3124\ (1,2) \\
\endmatrix$$
\midspace{.5\jot}
\caption{Table 1: Class of the permutation 1234 
in $S_4$ }

$$\spreadmatrixlines{.5\jot} 
\matrix  &4213 \ (3,4) &4132\ (4,5)  &4321\ (6,6) \\
         &2413 \ (2,4) &1432\ (5,5)  &3421\ (5,3) \\
         &2143 \ (4,4) &1342\ (3,2)  &3241\ (4,3) \\        
         &2134 \ (1,1) &1324\ (2,2)  &3214\ (3,3) \\
\endmatrix$$
\midspace{.5\jot}
\caption{Table 2: Class of the permutation 2134 
in $S_4$ }

We can also count permutations according to residue classes of 
major index with arbitrary modulus, by fixing the major index of the inverse 
modulo $(n-1).$ In fact,  Theorem 3.1 is the special case 
$k=n, \ell=n-1$ of the result below. 

\proclaim{Proposition 3.7}  Fix $2\leq k\leq n,$ $2\leq \ell$ 
a divisor of $(n-1),$ 
 and nonnegative integers $i\leq k-1,$ $j\leq n-2.$  Write $n=qk+r$ for 
$0\leq r\leq k-1.$  Fix an integer $s,$ $1\leq s\leq q.$
Then $s\cdot {n-1\over \ell} \cdot(n-2)!$ is 
the number of permutations $\tau$ in $S_n$ satisfying
  $$\maj(\tau)\equiv i\mod k, \ \ \maj(\tau^{-1})\equiv j\mod \ell,$$
and 
$$1\leq \tau^{-1}(n)\leq sk$$
\endproclaim

\demo{Proof} As in the proof of Proposition 3.5, 
we consider again the equivalence class of a permutation 
$\sigma$ in $S_{n-1}$ which fixes $(n-1).$  There are 
as usual $(n-2)!$ such classes.  Recall the rectangular 
representation of a class, with $n$ rows and $(n-1)$ columns. 
The permutation $\tau_{i,j}$ in row $i$ and column $j$ has the 
property that $n$ is in position $i,$ and when $n$ is erased, 
the resulting permutation in $S_{n-1}$ is the unique 
rearrangement of $\sigma$ with $(n-1)$ in position 
$(n-j+1),$ i.e., it equals $\sigma\ga_{n-1}^j.$

 Recall that 
the major index of the inverse is constant modulo $(n-1)$ on 
each column.  It is therefore constant modulo the divisor $\ell$ on 
each column.  Since the first row of our class is the set of 
circular rearrangements of $\sigma\in S_{n-1},$ with $n$ appended in position 
$n,$ we know that the major index of the inverse permutation in column $t$ 
is $t-1+\maj(\sigma^{-1})$ for $t=1,\ldots , n-1.$   In particular, 
for each fixed row $i,$ and for each $J=1, \ldots, {n-1\over \ell},$ 
columns $(J-1)\ell+1, \ldots, Jk,$ form a 
complete residue class modulo $\ell$ for the major index of the inverse 
permutation. Consequently in each fixed row $i,$ there are exactly 
${n-1 \over \ell}$ permutations $\tau_{i,t}$ whose inverse has major 
index congruent to $j\mod \ell.$ 

   Now consider the $sk(n-1)$ permutations $\tau_{i,j}$ for $i\leq sk,$ 
$1\leq j\leq n-1,$ in the 
sub-rectangle formed by the first $sk$ rows. 
We invoke Lemma 2.2. In the $sk$ by $(n-1)$ sub-rectangle formed 
by the first $sk$ rows, for each fixed $J,$  the major indices of 
the permutations in column $J$ 
form a set $M_J$ consisting of $sk$ consecutive integers.  This set $M_J$ 
can clearly be partitioned (as in the preceding paragraph) 
into $q$ disjoint subsets of consecutive integers, 
each of length $k,$ each of which is a complete residue system modulo $k.$  
Hence in each column $J,$ the first $sk$ rows contain exactly $s$ permutations 
$\tau_{i,J},\ i\leq sk, $ whose major index is congruent to $i\mod k.$

  Now all the permutations in row $i$ have $n$ appearing in
 position $i.$  
 Hence in the first $sk$ rows of the class of $\sigma$ 
we  find exactly ${n-1\over \ell}$  permutations $\tau_{i,j}$ 
 with the required properties.\qed 
\enddemo

By putting $r=0$ and $s=q={n\over k}$ in the above proposition, 
we immediately deduce 
the following elegant generalisation of Theorem 3.1.

\proclaim{Theorem  3.8} Let $k$ be a divisor of $n,$ $k\neq 1.$ 

Fix $0\leq i \leq k-1,\ 0\leq j\leq \ell-1.$  Then ${n!\over k\cdot\ell}$ 
is the 
 number of permutations $\tau$ in $S_n$ satisfying
  $$\maj(\tau)\equiv i\mod k, \ \ \maj(\tau^{-1})\equiv j\mod \ell.$$
\endproclaim
\endcomment

Let $m_n(i\bck k; j\bck\ell)$ denote the number of permutations 
$\sigma\in S_n$ 
with $\maj(\sigma)\equiv i\ \mod k$ and 
$\maj(\sigma^{-1})\equiv j\ \mod \ell.$   

\proclaim{Theorem 3.1}Let $\ell$ be a divisor of $n-1,$ $\ell\neq 1,$ and 
let $k$ be a divisor of $n,$  $k\neq 1.$ 
Fix $0\leq i \leq k-1,\ 0\leq j\leq \ell-1.$  Then
$$m_n(i\bck k; j\bck\ell)= {n!\over k\cdot\ell} .$$
\endproclaim

\demo{Proof} Let $\sigma\in S_{n-1},$ and construct $\sigma_i,$ 
$i=1,\ldots,n$ in $S_n$ as in Lemma 2.2, by inserting $n$ in position $i.$ 
Since $\ell|(n-1),$ we have by Lemma 2.4 that for all $i,$ 
$$\maj(\sigma^{-1})\equiv \maj(\sigma_i^{-1})\mod \ell.$$
  By Lemma 2.2, 
since the set $\{\maj(\sigma_i):i=1,\ldots,n\}$ consists of $n$ consecutive 
integers, each congruence class modulo $k$ appears exactly ${n\over k}$ times. 
Hence we have $$m_n(i\bck k; j\bck\ell)
={n\over k}\cdot |\{\sigma\in S_{n-1}:\maj(\sigma^{-1})\equiv j\mod \ell\}|$$
and the result now follows from Proposition 2.5.   \qed
\enddemo

By 
examining Lemma 2.2 more closely, we obtain the 
following recurrence on $n$ for these numbers in the case when $k$ and 
$\ell$ are divisors of $(n-1).$

\proclaim{Proposition 3.2} Let $\ell, k$ be  divisors of $n-1,$ $\ell\neq 1,$ 
$k\neq 1.$    Then 
$$m_n(i\bck k; j\bck\ell)=(n-2)! {(n-1)^2\over k\cdot\ell} 
+m_{n-1}(i\bck k; j\bck\ell).$$
\endproclaim

\demo{Proof} Let $\sigma\in S_{n-1},$ and construct $\sigma_i,$ 
$i=1,\ldots,n$ in $S_n$ as in Lemma 2.2, by inserting $n$ in position $i.$ 
Since $\ell|(n-1),$ we have by Lemma 2.4 that for all $i,$ 
$$\maj(\sigma^{-1})\equiv \maj(\sigma_i^{-1})\mod \ell.$$

 Now let $k|(n-1).$ 
 By Lemma 2.2, the major indices of the first $(n-1)$ elements 
$\sigma_i,$ $i=1,\ldots, n-1,$ form a segment of $(n-1)$ consecutive 
integers, and hence the residue class $i$ modulo $k$ appears 
exactly ${n-1\over k}$ times among them.  Also note that $\maj(\sigma)
=\maj(\sigma_n).$ 

Hence we have $m_n(i\bck k; j\bck \ell)$
$$={n-1\over k} |\{\sigma\in S_{n-1}:\maj(\sigma^{-1})
\equiv j\mod \ell\}|$$
$$+ |\{\sigma\in S_{n-1}:\maj(\sigma^{-1})
\equiv j\mod \ell, \ \maj(\sigma)\equiv i\mod k\}|.$$

Collecting terms and using Proposition 2.5, we obtain 
$$m_n(i\bck k; j\bck \ell)= {n-1\over k} {(n-1)!\over \ell}
+m_{n-1}(i\bck k; j\bck \ell),$$
\noindent as required.
\enddemo
\comment
\demo{Proof}  We follow the proof of Proposition 3.7 through the first 
two paragraphs, arriving at the conclusion that in each fixed row $i$ 
($i\leq n$), there are exactly ${n-1\over \ell}$ permutations 
$\tau_{i,t}$ whose inverse has major index congruent to $j\mod \ell.$ 
Now consider the first $n-1$ rows of the array 
$\{\tau_{i,j}\}_{1\leq i\leq n, 1\leq j\leq n-1}.$ In each column the 
major indices form a consecutive set of $(n-1)$ integers, by Lemma 
2.2, and hence form a complete residue system modulo $(n-1).$  They 
can also be partitioned into ${n-1\over k}$ complete sets of 
residue systems modulo $k,$ for $k$ dividing $(n-1).$ 
Hence we have that $(n-2)! {n-1\over \ell} {n-1\over k}$ 
equals the number of permutations $\tau$ in $S_n$ satisfying 
$$ \tau^{-1}(n)\neq n, \maj(\tau)\equiv i\mod k, \ 
\maj(\tau^{-1})\equiv j\mod \ell.$$ 

Since $m_n(i\bck k; j\bck\ell)$ counts  the above permutations 
plus those permutations in $S_n$ which fix $n$ and satisfy the 
same conditions on the major indices, we 
obtain the required result. \qed
\enddemo
\endcomment
We note that while the above arguments are not symmetric in $k$ and $\ell,$ 
 the numbers $m_n(i\bck k; j\bck\ell)$ satisfy 
$$m_n(i\bck k; j\bck\ell)=m_n(j\bck\ell; i\bck k).\tag C$$
This  follows by applying  the involution $\tau\mapsto \tau^{-1}.$ 

For arbitrary choices of $k, \ell,$ these  numbers usually depend 
on the values of $i$ and $j.$   For example for $n=4,$ we have 
$m_4(0\bck 2; 0\bck 2)=8=m_4(1\bck 2; 1\bck 2),$ 
and $m_4(1\bck 2; 0\bck 2)=4=m_4(0\bck 2; 1\bck 2).$ 
When $k=\ell=3,$ we have 
$m_4(0\bck 3;0\bck 3)=4,\ m_4(0\bck 3; 1\bck 3)=2=m_4(0\bck 3; 2\bck 3);$
and $m_4(1\bck 3;1\bck 3)=3=m_4(1\bck 3; 2\bck 3).$ 
The other values follow by symmetry from (C).

Note that in view of Proposition 2.5, we know that, for fixed 
$\ell,$ the sum over $i=0,1,\ldots , k-1$ of the numbers 
$m_n(i\bck k; j\bck \ell)$ is ${n!\over \ell}.$

\comment 
For each $n$ and each choice 
of integers $k,\ell\leq n,$ let $M_n(k, \ell)$ be the $k$ by $\ell$ 
 matrix with rows indexed 
by $i=0,1,\ldots, k-1,$ and columns indexed by $j=0,1,\ldots, \ell-1.$ 
The $(i,j)$th entry is the number $m_n(i\bck k; j\bck \ell)$ of 
permutations in $S_n$ with $\maj$ congruent to $i\mod k$ and $\maj$ of 
the inverse congruent to $j\mod \ell.$  

From Proposition 3.2, we obtain the following relation when $n-1=qk.$ 
$$m_n(i\bck k; j\bck k)
=(n-2)! q^2 +m_{n-1}(i\bck k; j\bck k).$$
In particular, 
$$m_n(i\bck n-1; j\bck n-1)
=(n-2)!  +m_{n-1}(i\bck n-1; j\bck n-1).$$

These observations are summarised in matrix terms below.
Part (1) below follows  from Proposition 2.5,  Part (2) follows 
from (C), and Part (3) follows from the preceding parts.  Part (4) 
is a special case of the preceding proposition.

\proclaim{Proposition 3.3} 
\item{(1)} The matrix $M_n(k,\ell)$ has constant column sum  equal to
 ${n!\over \ell},$ and constant row sum equal to 
${n!\over k}.$ 
\item{(2)} The matrices $M_n(k,\ell)$ are symmetric about the principal 
diagonal.  
\item{(3)} The matrix $M_n(2,2)$ has equal 
diagonal elements.
\item{(4)} Let $k|(n-1),$ and let $U_k$ denote the $k$ by $k$ 
matrix all of whose entries equal 1.  Then 
$M_n(k,k)=(n-2)! ({n-1\over k})^2 U_k + M_{n-1}(k,k).$

\endproclaim
\endcomment 
\comment
\proclaim{Proposition 3.3} Let $\ell|(n-1).$ 
 Let $k|n-2, k\geq 2.$  Then
for all $0\leq i_0\leq k-1$ and $0\leq j\leq \ell-1,$
$$m_n(i_0\bck k;j\bck \ell)={n!\over k\cdot \ell}.$$
\endproclaim

\demo{Proof} Let $M_n(i\bck k;j\bck\ell)$ denote the set of permutations
$\tau\in S_n$ such that $\maj(\tau)\equiv i\mod k,$
$\maj(\tau^{-1})\equiv j\mod \ell.$   
Then $m_n(i\bck k;j\bck\ell)$ is the cardinality
 of $M_n(i\bck k;j\bck\ell).$

Consider $\tau\in M_{n-1}(i\bck k;j\bck\ell),$ and let $\tau_t$ denote 
the permutation in $S_n$ obtained by inserting $n$ into $\tau$ in
position $t,$ for $t=1,\ldots , n.$   Since $\ell|n-1,$ by Lemma 2.4
$\maj(\tau_t^{-1})=\maj(\tau^{-1})\equiv j\mod \ell$ for all $t.$

 Now assume $k|(n-2).$
By Lemma 2.2, the set $\{\maj(\tau_t):1\leq t\leq n-2\}$ consists of
$(n-2)$ consecutive integers, and therefore forms a complete
residue system modulo $(n-2).$   Since $k|(n-2),$ each residue class
modulo $k$ appears exactly $(n-2)/k$ times in this set.
Now consider $\tau_n$ and $\tau_{n-1}.$  Clearly $\tau_n$ has the
same major index as $\tau,$  congruent to $i\mod k.$
It is easy to see that $\maj(\tau_{n-1})-\maj(\tau)$ equals $(n-1)$ or
1, depending on whether or not $(n-2)$ is a descent of $\tau.$
  In either case, $\maj(\tau)$
is congruent to $(i+1)$ modulo $n-2,$ and hence is congruent to $(i+1)$
modulo $k.$

These arguments show that in the set $\{\maj(\tau_t):1\leq t\leq n\},$ 
the residue classes $i$ and  $i+1$ modulo $k$ each appear exactly $(n-2)/k+1$
times, while the remaining $(k-2)$ residue classes modulo $k$ 
appear exactly $(n-2)/k$ times.  

Hence for any fixed residue class $i_0$ modulo $k,$ we have 
that $m_n(i_0\bck k;j\bck \ell)$ equals 
$$({n-2\over k}+1)(m_{n-1}(i_0\bck k;j\bck\ell)+
m_{n-1}(i_0-1\bck k;j\bck\ell))
 +{n-2\over k}
\sum_{i\neq i_0, i\neq i_0-1}m_{n-1}(i\bck k;j\bck\ell),$$
\noindent for all $j.$

Collecting terms,  the right-hand side reduces to 
$$m_{n-1}(i_0\bck k;j\bck \ell)+m_{n-1}(i_0-1\bck k;j\bck \ell) 
+{n-2\over k}\sum_{i=0}^{k-1} m_{n-1}(i\bck k;j\bck \ell).$$

The sum $\sum_{i=0}^{k-1} m_{n-1}(i\bck k;j\bck\ell)$ is simply the number of 
permutations in $S_{n-1}$ with major index of the inverse 
congruent to $j\mod\ell,$ and hence by Proposition 2.6 it equals 
$(n-1)!/\ell.$  Also by observation (C) and Proposition 3.1 (1), 
the first 
two terms are each equal to ${(n-1)!\over k\cdot\ell},$ since 
$k|n-2$ and $\ell|n-1.$  Hence we have 

$$m_n(i_0\bck k;j\bck\ell)=2 {(n-1)!\over k\cdot\ell} 
+{n-2\over k}{(n-1)!\over \ell}
={n!\over k\cdot \ell},$$
\noindent as required. \qed
\enddemo
\endcomment
\comment
From the preceding results we can conclude
\proclaim{Theorem 3.4} Let $k<\ell\leq n.$  Then for all $i,j,$ we have 
$$m_n(i\bck k;j\bck \ell)={n!\over k\cdot\ell}$$
\noindent whenever $\ell|(n-1)$ and $k|n$ or $k|(n-2).$ 
\endproclaim
\endcomment

\head {4. Some bijections}\endhead

In this section we present bijective proofs for some of the 
results derived in Sections 3 and 2.
Recall that this paper was originally motivated by the algebraic 
discovery of the formula (A).  We now give a bijective proof of (A), 
which is the special case $k=n,\ell=n-1$ of Theorem 3.1.

\proclaim{Proposition 4.1} (Bijection for the case $k=n,\ell=n-1$ 
of Theorem 3.1.)  Fix integers  $0\leq i\leq n-1, 0\leq j\leq n-2.$ 
Then the number of permutations $\sigma$ in $S_n$ such that 
$\maj(\sigma)\equiv i\mod n$ and $ \maj(\sigma^{-1})\equiv j\mod (n-1),$ 
equals $(n-2)!$
\endproclaim

\demo{Proof} First note that $(n-2)!$ counts the number of permutations 
in $S_{n-1}$ having $(n-1)$ as a fixed point.  Let $A_{n-1}$ be this set of 
permutations, and let $B_n$ be the subset of $S_n$ with major indices as 
prescribed in the statement of the 
theorem.  Given $\sigma\in A_{n-1},$ by Lemma 2.1 there is a unique circular 
rearrangement $\sigma'$ in $S_{n-1}$ whose inverse has 
 major index congruent to $j \mod (n-1).$  Lemma 2.2 then 
shows that, for each $i=0,1,\ldots, n-1,$  
there is a unique position in $\sigma'$ in which to insert $n,$ in order to  
obtain a permutation $\sigma''\in S_n$ such that 
$\maj(\sigma'')\equiv i \mod n.$   By Lemma 2.4, the passage from $\sigma'$ 
to $\sigma''$ does not change the major index of the inverses modulo 
$(n-1),$ and thus 
$\maj({\sigma''}^{-1})=\maj({\sigma'}^{-1})\equiv j \mod (n-1).$ 
Hence  $\sigma\mapsto \sigma''$ gives a well-defined map from 
$A_{n-1}$ to $B_n.$ 
To see that this is a bijection, given $\sigma''\in B_n,$ erase the $n$ to 
obtain $\sigma'\in S_{n-1},$ and let $\sigma$ be the unique 
circular rearrangement of $\sigma'$ such that $\sigma(n-1)=n-1.$ 
Then $\sigma\in A_{n-1},$ and clearly the map is a bijection. \qed

\enddemo

\proclaim{Example 4.1.1} Let $n=6, i=2, j=3.$ Take $\sigma=21345\in A_5.$ 
Note that $\maj(\sigma^{-1})=1.$ The unique circular rearrangement 
whose inverse has major index equal to $3(\equiv 3 \mod 5)$ is 
$\sigma'=34521.$ Now $\maj(\sigma')=7,$ (descents in positions 3 and 4).
Now use (the proof of) Lemma 2.2. 
To obtain a permutation with major index 8 ($\equiv 2\mod 6$), insert 
$6$ into position 5 (immediately after the right-most descent).
This gives $\sigma''=345261\in B_6.$
\endproclaim
 
The remainder of this section is devoted to giving a constructive  proof 
of Proposition 2.5.    A bijection for the case $k=n$ was given in 
\cite{1}, using Lemma 2.1.  We do not know of a bijection for arbitrary 
$k,$ but a bijection for the case $k=n-1$ is given in the proof 
which follows.

\proclaim{Proposition 4.2} (Bijection for the case $k=n-1$ 
of Proposition 2.5.)
 Fix an integer $0\leq j\leq n-2.$  
The number of permutations in $S_n$ with major index congruent to 
$j\mod (n-1)$ is $n(n-2)!=n!/(n-1).$
\endproclaim

\demo{Proof} Let $B_n$ denote the set $\{\sigma\in S_n: \maj(\sigma^{-1})
\equiv j\mod (n-1)\}.$  It suffices to show that this set has 
cardinality $n(n-2)!$   Let $C_n$ denote the 
set of permutations $\tau\in S_n$ such that, when $n$ is erased, 
$(n-1)$ is a fixed point of the resulting permutation $\tau'$ in $S_{n-1}.$ 
Observe that $C_n$ has cardinality $n(n-2)!,$ since the 
number of permutations in $S_{n-1}$ which fix $(n-1)$ is $(n-2)!,$ and 
there are $n$ positions in which $n$ can be inserted.

We describe a bijection between $C_n$ and $B_n.$  If $\tau\in C_n,$ 
let $\tau'$ be the permutation in $S_{n-1}$ obtained by erasing $n.$ 
By definition of $C_n,$ $\tau'(n-1)=n-1.$  By Lemma 2.1, there is a 
unique circular rearrangement $\tau''\in S_{n-1}$ of $\tau'$ such 
that the major index of the inverse of $\tau''$ 
is congruent to $j\mod (n-1).$  Now construct $\tilde{\tau}\in S_n$ by 
inserting  $n$ into $\tau''$ in the same position that it occupied 
in $\tau,$ i.e., ${\tilde{\tau}}^{-1}(n)=\tau^{-1}(n).$  By Lemma 2.4, 
$\maj({\tilde{\tau}}^{-1})=\maj({\tau''}^{-1})\equiv j\mod (n-1).$ 
Hence we have 
a map $\tau\mapsto \tilde{\tau}\in B_n.$ It is easy to see that this 
construction can be reversed exactly as in the proof of Proposition 4.1, and 
hence we have the desired bijection. \qed
\enddemo

\proclaim{Example 4.2.1}Let $n=5, j=2.$ Take $\tau=32154.$ 
Then $\tau$ belongs to the set $ C_5$ of the preceding proof.  
Erasing 5 yields $\tau'=3214,$ whose inverse major index is 3. 
The third cyclic rearrangement $\tau''=4321$ then has inverse major
 index 6 $\equiv 2\mod 4,$ and $\tau\mapsto \tilde{\tau}=43251.$

\endproclaim

Now we examine Klyachko's Lemma 2.1 more closely.   We obtain
the following  result, which specialises, in the case $k=n, $ to 
 Proposition 1.3.

\proclaim{Lemma 4.3} Fix integers $1\leq k\leq n,$ $0\leq j\leq k-1$ 
and $1\leq a \leq n-k+1.$
\item{(1)}
Then $$(n-1)!=|\{\sigma\in S_n:\maj(\sigma^{-1})\equiv j\mod k {\text\ and\ } 
n-a-k+2\leq \sigma^{-1}(n)\leq  n-a+1.\}|$$
\item{(2)} Let $n=qk+r,$ $0\leq r\leq k-1.$
 Fix an integer $s$ between 1 and $q.$ 
Then $$s(n-1)!=\{\sigma\in S_n: maj(\sigma^{-1})
\equiv j\mod k {\text \ and\ } \sigma^{-1}(n)\in [n-sk+1,n]\}.$$

\endproclaim

\demo{Proof}  Let $A_n$ denote the set of permutations in $S_n$ which fix $n,$ 
and let $B_n$ denote the subset of $S_n$ subject to the conditions 
in the statement of Part (1).
Let $\tau\in A_n.$  Consider the circular class of $\tau$ consisting 
 of the set $\{\tau, \tau\ga, \ldots, \tau\ga^{n-1}\}.$  The proof of 
Lemma 2.1 shows that because  $\tau(n)=n,$ 
we have the exact equality $\maj(\tau\ga^i)=\maj(\tau)+i,$ for 
$0\le i\le n-1.$  In particular, for any $1\leq k\leq n,$ the first $k$ 
circular rearrangements $\tau\ga^i, 0\leq i\leq k-1,$ have the property that 
the major indices of their inverses form a complete 
residue system modulo $k.$ More generally, this observation holds 
for any $k$ consecutive circular rearrangements $\tau\ga^i, a\leq i\leq a+k-1,$
where $a$ is any fixed integer $1\leq a \leq n-k+1.$ 

Hence  for every $\tau\in A_n,$ there is a unique $i, a \leq i\le a+k-1$ such 
$\sigma=\tau\ga^i$ has $\maj(\sigma^{-1})\equiv j \mod k.$ 
Since $n$ is in position $n-i$ in $\tau\ga^i,$ clearly 
$n-a-k+2\leq \sigma^{-1}(n)\leq n-a+1.$
Thus $\tau\mapsto\sigma $ gives a well-defined map from $A_n$ to $B_n.$
Conversely given $\sigma\in B_n,$ with 
$\sigma^{-1}(n)=n-i+1,\ a\leq i\leq a+k-1,$ 
let $\tau\in S_n$ be defined by $\tau\ga^i=\sigma.$ Then 
clearly $\tau(n)=n,$  and $\tau\in A_n.$  This shows that our map is a 
bijection, and (1) is proved.

For (2), again we start with the set $A_n$ of the $(n-1)!$ 
permutations in $S_n$ which fix $n.$ Let $\tau\in A_n.$ Then as in 
the preceding proof, for $i=0,1,\ldots, sk-1,$ the first $sk$ 
circular rearrangements $\tau\ga^i$ have $n$ in position $(n-i),$ and 
$\maj((\tau\ga^i)^{-1})=\maj(\tau)+i.$  In particular, for each 
$J=1,\ldots, s,$ the major index of the inverse permutations in the 
subset $\{\tau\ga^{(J-1)k+i}:0\leq i\leq k-1\}$ is a complete residue 
system modulo $k.$ Hence the first $sk$ rearrangements contain 
exactly $s$ permutations with inverse major index congruent to $j\mod k.$
This establishes (2).  \qed

\enddemo

We are now ready to give a constructive proof of an equivalent restatement of 
Proposition 2.5, by looking at the circular classes of permutations 
$\tau\in S_n$ which fix $n.$  Note that the statement of Proposition 4.4
 ( or Proposition 2.5) is invariant  with respect to taking inverses, 
i.e., it says that ${n!\over k}$ 
is also the number of permutations in $S_n$ with constant 
major index modulo $k.$  Our constructive proof, however, works only 
for the inverse permutations.

\proclaim{Proposition 4.4} (Constructive proof)
$${n!\over k}=|\{\sigma\in S_n: \maj(\sigma^{-1})\equiv j\mod k\}|.$$
\endproclaim

\demo{Proof}    We proceed 
inductively.  We assume $k\leq n-1,$ since the case $k=n$ was dealt with 
in Proposition 1.3.  It is easy to verify directly that the statement holds 
for $n=3.$ Assume we have constructed the permutations in $S_{n-1}$ 
with inverse major index congruent to $j\mod k.$ Note that this means 
we can identify these permutations in the subset $A_n$ of $S_n.$ 
Let $\tau\in A_n.$  We show how to pick out the permutations in 
the circular class of $\tau$ with inverse major index congruent to $j\mod k.$ 
Let $n=qk+r.$  Taking $s=q$ in Lemma 4.3 (2), the proof shows how to 
pick out the $q$ permutations in the first $qk$ circular rearrangements 
$\tau\ga^i, 0\leq i\leq qk-1.$ Now consider the remaining $r$ (recall $r<k$)  
rearrangements $\tau\ga^i, qk\leq i\leq qk+r-1.$  These will contain 
a (necessarily unique) permutation with inverse major index 
congruent to $j\mod k,$ 
iff $\maj(\tau^{-1})\equiv j-i\mod k,$ for $qk\leq i\leq qk+r-1,$
i.e., iff $\maj(\tau^{-1})\equiv j-t\mod k,$ for $t=0,\ldots ,r-1.$ 
By induction hypothesis for each $t=0,\ldots, r-1,$ there are exactly 
$(n-1)!/k$ such permutations in $A_n.$ Hence there are $r (n-1)!/k$ 
 permutations in $A_n$ whose circular class is such that, among the last $r$ 
rearrangements, there is  a permutation with inverse $\maj$ congruent to 
$j\mod k.$ 

We have thus accounted for exactly $q(n-1)! +r(n-1)!/k=n!/k$ 
permutations $\sigma\in S_n$ 
with $\maj(\sigma^{-1})\equiv j \mod k.$  \qed
\enddemo



\refstyle{C}

\enddocument